\documentclass[a4paper, 12pt]{article}

\usepackage[dvips]{graphics,color}

\usepackage{amsmath, amssymb}
\usepackage{verbatim}
\usepackage[british]{babel}
\usepackage{graphicx}

\usepackage{setspace}

\newtheorem {lemma}{Lemma}
\newtheorem {corollary}{Corollary}
\newtheorem {prop}{Proposition}
\newtheorem {theorem}{Theorem}
\newtheorem{rem}{Remark}

\newcommand{\R}{{\mathbb{R}}}
\renewcommand{\P}{\mathsf{P}}
\newcommand{\Z}{\mathbb{Z}}
\newcommand{\N}{\mathbb{N}}
\newcommand{\eps}{\varepsilon}
\newcommand{\E}{\mathsf{E}}
\newcommand{\G}{\mathsf{L}}

\date{}

\begin{document}

\title{Long term behaviour  of two interacting birth-and-death processes}

\author{
Mikhail Menshikov\footnote{Department of Mathematical Sciences,
 Durham University, UK.
 Email address: mikhail.menshikov@durham.ac.uk
}\\
{\small  Durham University}
\and 
Vadim Shcherbakov\footnote{Department of Mathematics, Royal Holloway,  University of London, UK.
 Email address: vadim.shcherbakov@rhul.ac.uk
}\\
{\small  Royal Holloway,  University of London}
}

\maketitle

\begin{abstract}
{\small 
In this paper we study  the long term evolution of a continuous time  Markov chain  formed by  
two interacting birth-and-death processes.  The interaction between the processes 
is modelled by  transition  rates which are functions with suitable monotonicity properties. 
This is in  line with   the approach 
proposed by    Gauss G.F. and Kolmogorov   A.N. for modelling interaction between species in ecology.
We obtain conditions for transience/recurrence of the Markov chain and describe in detail its asymptotic 
behaviour  in special   transient   cases. In particular, we find that in some of these  cases  the Markov chain escapes to infinity
in an unusual way, and the corresponding trajectories   can be rather precisely  described.  
 }
\end{abstract}


\section{Introduction}
\label{growth}

A birth-and-death process on  $\Z_{+}=\{0,1, 2, \ldots\}$
   is a continuous time Markov chain (CTMC) that evolves as follows. 
Given a current  state  $k$ it jumps either to $k+1$,  or to $k-1$ (if  $k>0$) at certain state dependent rates.
The long term behaviour of a  birth-and-death process   is well known. 
Namely, given a set  of transition  rates   one can, in principle, determine  whether the corresponding 
birth-and-death process   is recurrent/positive recurrent, or  transient/explosive,  and compute various characteristics of the process. 
These results can be found in many books  (e.g., see \cite{Feller}, \cite{Karlin0} and  \cite{Liggett}).
The long term behaviour of multivariate Markov processes with similar  dynamics  is less known. 

In this paper we study  the long term behaviour of 
 CTMC  $\xi(t)=(\xi_1(t), \xi_2(t))\in \Z_{+}^2,$ 
 evolving as follows.
Given $\xi(t)=(x, y)\in \Z_{+}^{2}$  the Markov chain jumps to $(x+1, y)$ and  to  $(x, y+1)$
 at rates $F(x)G(y)$ and  $F(y)G(x)$ respectively, 
 where $F$ and $G$ are positive functions on  $\R_{+}=[0, \infty)$.
Also, the Markov chain   jumps   from $(x, y)$   to  $(x-1, y)$ at the constant rate of $1$, provided $x>0$, and 
it jumps  to $(x, y-1)$ at  the same constant rate of $1$, provided that $y>0$. 

The Markov  chain is a two-dimensional analogue of integer valued  birth-and-death processes, and 
  can be interpreted in terms of  two interacting birth-and-death processes. 
The construction of the  birth rates  allows to model various types of both individual dynamics and interaction between the Markov chain components. 
Function $F$ determines, in terms of statistical physics,   the free dynamics of a component.
Interaction between components is  modelled by choosing  an appropriate function  $G$.
If, say,   $G\equiv 1$, then  $\xi_1(t)$ and $\xi_2(t)$ are  independent identically distributed birth-and-death processes.
Given $F$,  one can  choose  a  decreasing $G$ in order to model a  competitive interaction.
If $G$ is  increasing, then  a component's growth is accelerated by  its  neighbour.

Recall that a birth-and-death process on $\Z_{+}$   is a classic probabilistic  model for the size of a population.
Therefore  CTMC $\xi(t)$ can be regarded   as a stochastic model for  two interacting populations.
The model is  related  to  stochastic population models 
 formulated  in terms of  two interacting birth-and-death processes  (e.g., see \cite{Anderson}, \cite{Klebaner}, \cite{Becker}, 
\cite{Reuter}, \cite{Ridler} and references therein). In these 
models, which are stochastic versions of   the famous  Lotke-Volterra model, a pair of birth-and-death processes typically evolves as follows. 
 Given a current state $(x, y)$ of the components,  the  individual transition rates  are linear in $x$ and $y$, while  interaction terms, included usually in death rates (i.e. competitive interaction), are  proportional to  $xy$. 
Our model is  in the spirit of the more general  approach proposed by Gauss G.F. (\cite{Gauss})
and Kolmogorov A.N. (\cite{Kolm72}) for modelling interactions in ecology.
Although they considered  deterministic population models, the idea is rather general. 
 According to  this approach, the interaction between species should be  modelled by   transition rates specified 
 by  general   functions with suitable (suggested by a motivating application)  monotonicity properties. 
A brief,  but  informative presentation  of these ideas is given in  \cite{Sigmund}, where further  references 
can be found. 
In our model the  interaction is built into the birth rates, though the model can be generalised by allowing  for non-constant death rates. 
We do not  explore further the  relationship of the Markov chain  with stochastic population models
and  focus on its  long term behaviour which is of interest from a mathematical point of view.

 It should be also noted that our Markov chain  is a  particular example
 of non-homogeneous random walks. 
The long term behaviour of non-homogeneous random walks  is much less studied (e.g., see  \cite{MPW}  and references therein)
 in contrast to homogeneous random walks  in domains with boundaries (e.g., see \cite{FMM} and references therein). 

We systematically apply the Lyapunov function  approach in our proofs. This approach  is well known and widely used 
for determining whether a Markov process is recurrent or transient (e.g., see \cite{FMM}, \cite{MPW} and references therein). 
 In Theorem~\ref{T1} we establish whether  the Markov chain is  transient or recurrent  under fairly general assumptions on functions $F$ and $G$.
Though the asymptotic behaviour of the Markov chain in this theorem can be guessed  from approximate sketches of the vector field of 
 its  mean infinitesimal jumps (see Figures~\ref{Fig1} and~\ref{Fig2}),  the Lyapunov function approach helps
 to formalize these  intuitive ideas. 
In Theorems~\ref{T2} and~\ref{T3} we obtain a  more detailed description of the long term behaviour of the Markov chain in some 
 transient cases. 
It should be noted that the Lyapunov function method  is also a powerful tool 
 for detecting phenomena that  might not be  immediately visible and are more  refined than just
 recurrence/transience.  Theorem 3  below provides an example of such a phenomenon. 
 In particular we show that in a  transient case specified by polynomial functions $F$ and $G$ 
  the Markov chain with probability one escapes to infinity in the following way. 
Namely,  the Markov chain is eventually absorbed  to either 
a horizontal strip $\{(x,y): y\leq k\}$, or a vertical strip $\{(x,y): y\leq k\}$, where  $k$ is explicitly computable.
Moreover, being eventually adsorbed by the  horizontal (vertical)  strip,   the Markov chain visits every line $y=i,\, i=0, \ldots, k$
( $x=i,\, i=0, \ldots, k$) infinitely often.

\section{Results}

Let  $(\Omega, {\cal F}, \P)$ be a probability space on which the Markov chain is defined.
Denote by $\E$ the expectation with respect to probability  measure $\P$.
Recall that the embedded Markov chain, corresponding to a CTMC, is a discrete time Markov chain (DTMC) 
with the same state space, and that makes the same jumps as the CTMC with probabilities proportional to the corresponding jump rates. 
Let  $\zeta(t)=(\zeta_1(t), \zeta_2(t))\in \Z_{+}^2$ be the DTMC  corresponding to the CTMC $\xi(t)$.  
Note that  we  use the same symbol $t$ for discrete time.
Given a real valued function $f$ on $\Z_{+}^2$  denote
\begin{equation}
\label{mf}
m_f(x, y, t)=\E(f(\zeta_1(t), \zeta_2(t))|\zeta(0)=(x, y)) -f(x, y),\quad (x, y)\in Z_{+}^2, \, t\in \Z_{+}.
\end{equation}
It is easy to see that
\begin{equation}
\label{mf1}
m_f(x, y, 1)=\frac{\G f(x, y)}{\gamma(x, y)},
\end{equation}
where 
\begin{align}
\label{Gener}
\G f(x,y)&=(f(x+1, y)-f(x,y))F(x)G(y)+(f(x-1, y)-f(x,y))1_{\{x>0\}}\\
&+(f(x, y+1)-f(x,y))F(y)G(x)+(f(x, y-1)-f(x,y))1_{\{y>0\}}\nonumber
\end{align}
is the generator of CTMC $\xi(t)$, and
\begin{equation}
\label{gamma}
\gamma(x, y)=F(x)G(y)+F(y)G(x)+1_{\{x>0\}}+1_{\{y>0\}},
\end{equation}
is the total intensity of jumps of CTMC $\xi(t)$.
In the last two equations and in what follows, $1_{A}$ denotes the indicator function of a set $A$. 
Note   that $\gamma(x, y)=\gamma(y, x)$.

   Recall that   a real valued function $g$  is called non-decreasing (non-increasing)
 on a set $A\subseteq \R$, if  $g(x)\leq g(y)$ ($g(x)\geq g(y)$)  for all $x, y\in A$,
such that $x\leq y$. 
Finally, throughout the text  we denote by $C_i,\, i=1,2,...$, or, just $C$, various constants, whose exact values 
are immaterial.

We are ready now to formulate the findings of our paper. We start with the classification of the long term behaviour of the Markov chain 
under fairly general assumptions on functions $F$ and $G$. 
\begin{theorem}
\label{T1}
Let  functions $F$ and $G$ be  positive. 
\begin{itemize}
\item[ 1)]  Let function  $F$ be  non-increasing and  $\lim_{x\to \infty}F(x)=0$.

a) If  one of the following two assumptions holds
\begin{itemize}
\item 
function $G$ is non-increasing and  $\lim_{x\to \infty}G(x)=0$, 
\item 
function  $G$ is non-decreasing, $\lim_{x\to \infty}G(x)=\infty$ and 
 $\lim_{x\to \infty}F(x-1)G(x)=0$,   
\end{itemize}
then CTMC $\xi(t)$ is   positive recurrent.

b)  If  function   $G$ is non-decreasing and $\lim_{x\to \infty}F(x)G(x)=\infty$,   then  CTMC $\xi(t)$  is  transient.

\item[ 2)] If function $F$ is non-decreasing, $\lim_{x\to \infty}F(x)=\infty$ and 
 one of the   following two  assumptions holds  
\begin{itemize}
\item function $G$ is non-decreasing and $\lim_{x\to \infty}G(x)=\infty$,
\item   function $G$ is non-increasing, $\lim_{x\to \infty}G(x)=0$,  and \newline 
either 
 $\lim_{x\to \infty}F(x)G(x)=0$, or $\lim_{x\to \infty}F(x)G(x)=\infty$, 
\end{itemize}
then  CTMC $\xi(t)$ is transient.

\end{itemize}
\end{theorem}

\begin{figure}[htbp]
\centering
\begin{tabular}{cc}
   \resizebox{0.4\textwidth}{!}{\includegraphics{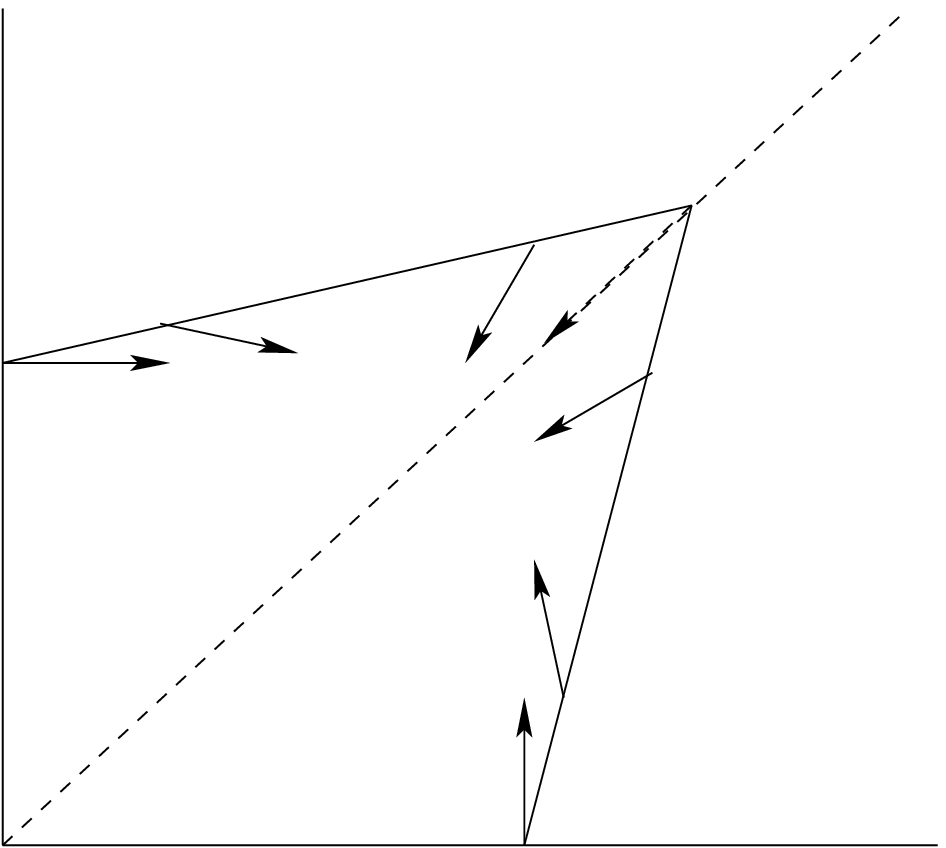}}&
     \resizebox{0.4\textwidth}{!}{\includegraphics{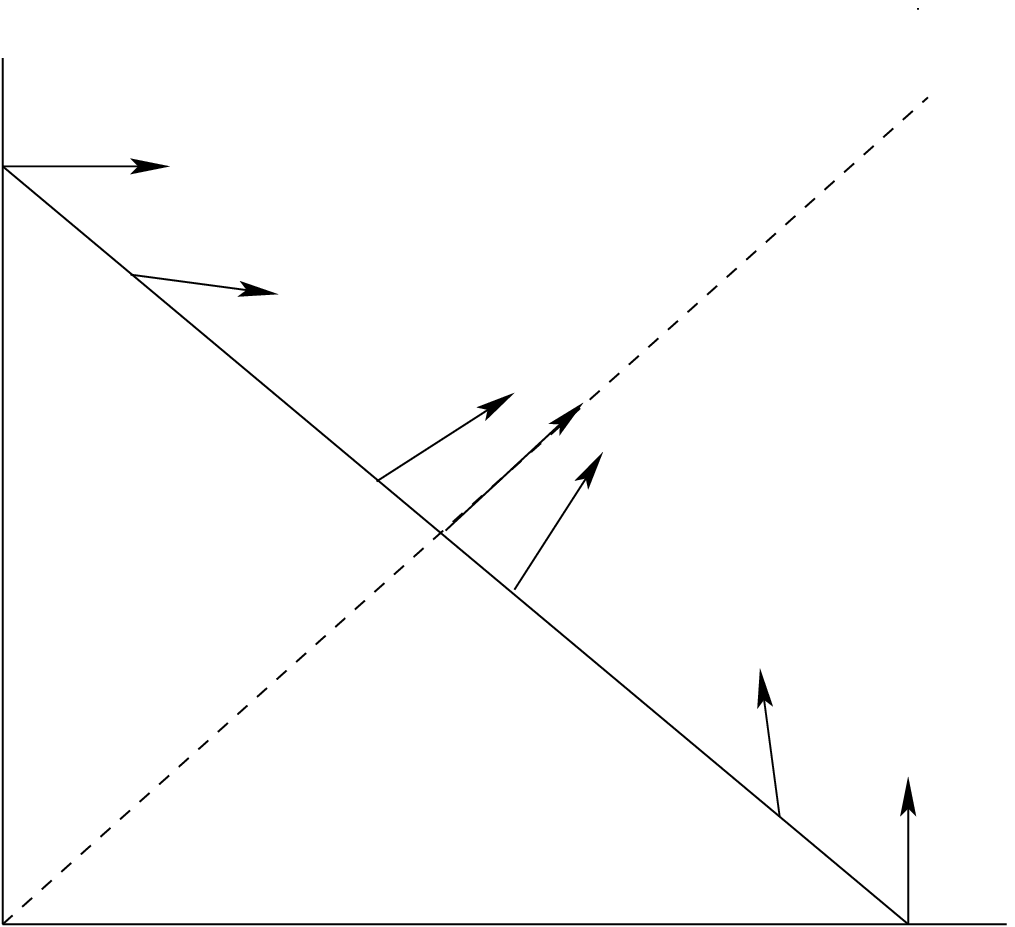}} \\
\end{tabular}
\caption{{\footnotesize 
 $\lim_{x\to \infty} F(x)=0$,\,  $\lim_{x\to \infty} G(x)=\infty$;
Left:  $\lim_{x\to \infty}F(x-1)G(x)=0$;
Right:  $\lim_{x\to \infty}F(x)G(x)=\infty$.
}}
\label{Fig1}
\end{figure}

\begin{figure}[ht]
\begin{center}
\begin{tabular}{cc}
   \includegraphics[width=2.7in, height=2.2in]{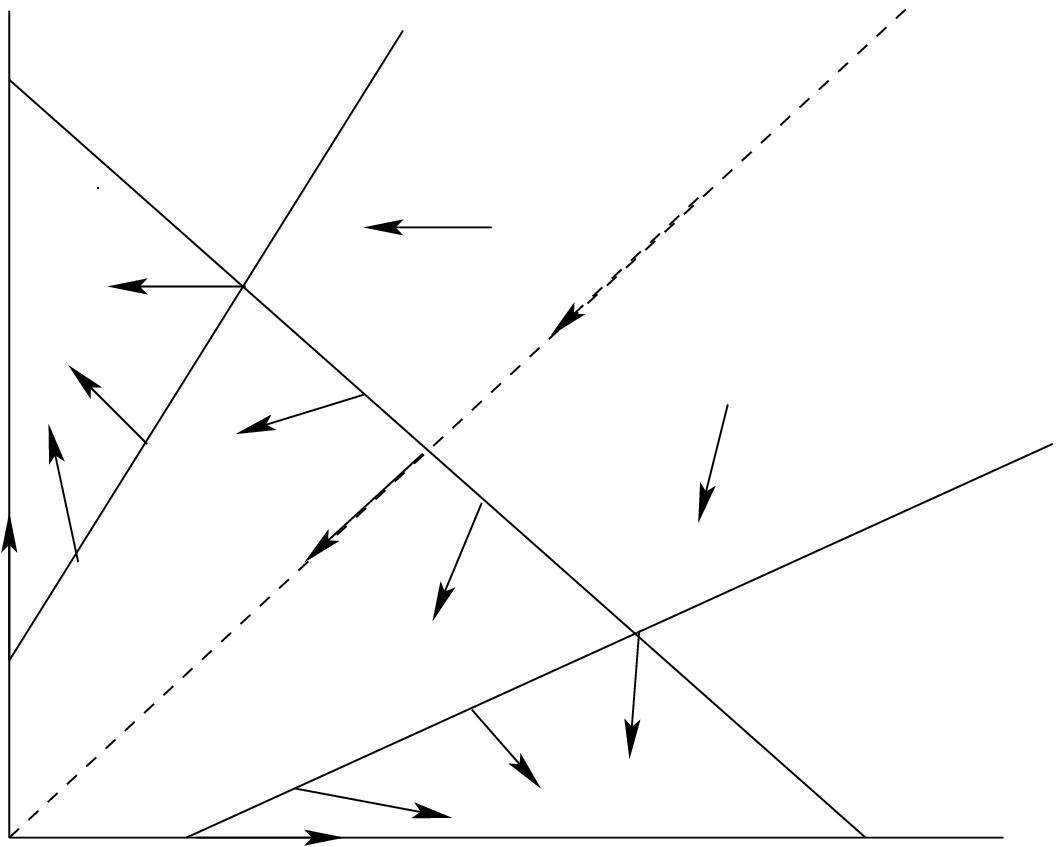}&
     \includegraphics[width=2.7in, height=2.2in]{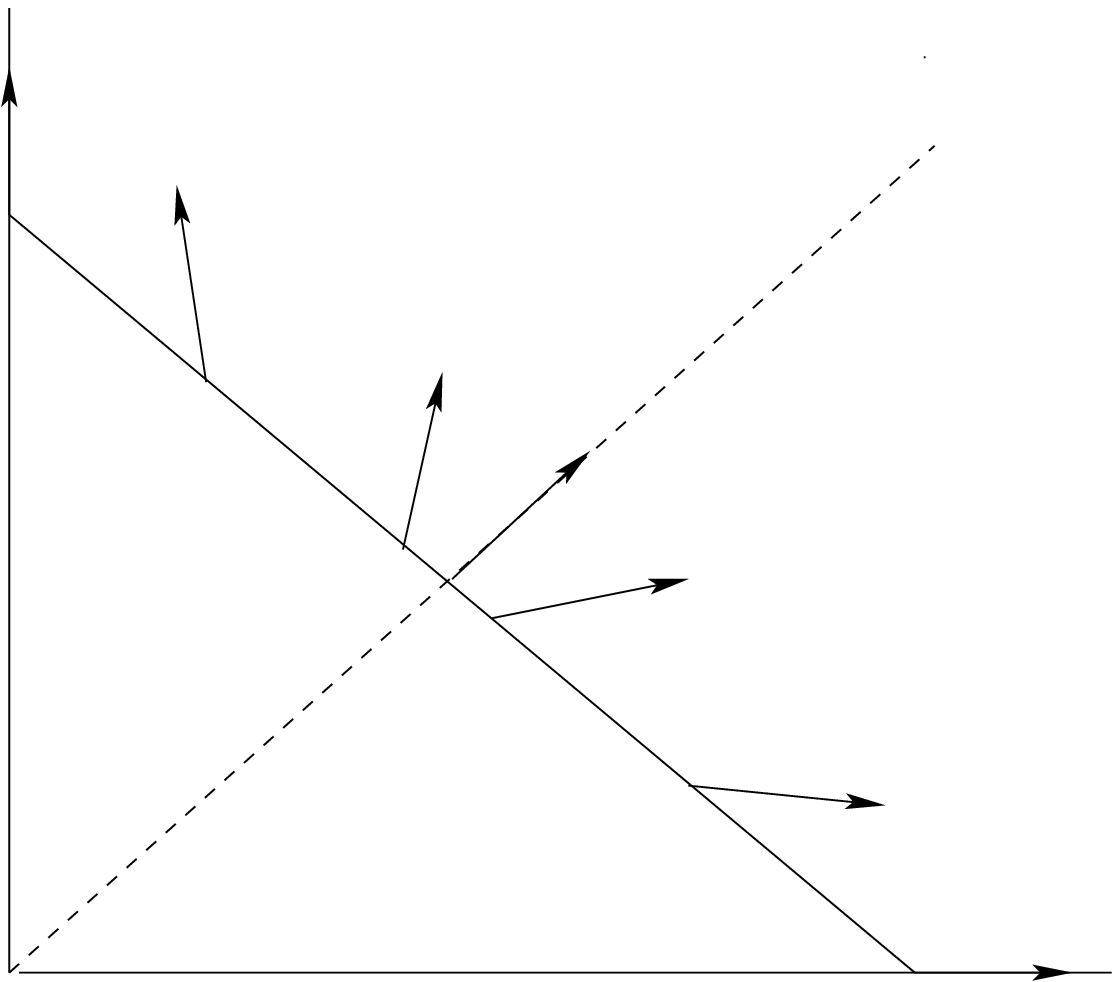} \\
\end{tabular}
\caption{{\footnotesize $\lim_{x\to \infty} F(x)=\infty$,\,  $\lim_{x\to \infty} G(x)=0$;
Left:  $\lim_{x\to \infty}F(x)G(x)=0$;
Right:  $\lim_{x\to \infty}F(x)G(x)=\infty$.
}}
\label{Fig2}
\end{center}
\end{figure}

\begin{rem}
\label{R1}
{\rm 
It is easy to see that Theorem~\ref{T1} describes the long term evolution of the Markov chain 
in  six  different cases.  Firstly,  if both $F(x)$ and $G(x)$  are non-increasing and have limit $0$ 
at infinity (and, hence,  $F(x)G(x)\to 0$ as $x\to \infty$), then  the Markov chain is  positive recurrent.
Secondly, if  both $F(x)$ and $G(x)$ increase 
to infinity (and, hence,  $F(x)G(x)\to \infty$ as $x\to \infty$), then   the Markov chain  is transient.  
  Approximate sketches of 
 a vector field of mean infinitesimal jumps of the Markov chain in other four cases are shown in Figure~\ref{Fig1} and Figure~\ref{Fig2}.
}
 \end{rem} 

\begin{rem}
\label{R2}
{\rm 
It should be noted  that  assumptions of the theorem
 are mostly  motivated  by the case of polynomial functions, e.g.  $F(x)=(x+1)^{\alpha}$,
 $\alpha\in \R$, and  $G(x)=(x+1)^{\beta}$, $\beta\in \R$.
Some of these assumptions can be slightly weakened without changing the proof. For example, in Part 2) the infinite limit of the product 
$FG$ at infinity in the case of non-increasing $G$ can be replaced by a sufficiently large limit (at least $2$). Such generalizations are not 
of much  interest.
 Also, some of these assumptions can be weakened provided that an additional information is available about functions $F$ and $G$
(e.g. see Remark~\ref{R5} in Appendix).}
\end{rem}

\begin{rem}
\label{R3}
{\rm 
Let us also  discuss  assumption ({\bf A1}): $\lim_{x\to \infty}F(x-1)G(x)=0$ in Part 1)a) of the theorem.
Ideally, we would like to replace it by  the following  assumption ({\bf A2}): $\lim_{x\to \infty}F(x)G(x)=0$.
Assumption ({\bf A1}) is violated, for example,  by functions  $F(x)=e^{-x^2}$ and $G(x)=e^{x^2}/x$.
Note that   assumptions ({\bf A1}) and ({\bf A2}) are equivalent in many cases. Moreover, in many cases 
these assumptions  are equivalent to the following  stronger assumption 
  ({\bf A3}): $\lim_{x\to \infty}F(\gamma x)G(x)=0$, where $\gamma \in (0, 1)$.
For example, this is the case if  $F(x)$ is a regularly  varying function of index $\alpha<0$. 
Equivalence can take place for  a non-regular varying $F$  as well, for example, 
if $F(x)=e^{-\alpha x}$ and $G(x)\leq e^{\beta x}$, where  $\alpha, \beta>0$ and $\alpha>\beta$.
}
\end{rem}

\begin{rem}
 \label{R4}
{\rm 
It should be noted   that there is 
 a certain  phase transition in the long term  behaviour of the Markov chain in the case of non-increasing and vanishing at infinity $F$. 
Indeed, if $G$ is also non-increasing with zero limit at infinity, then  the Markov chain is  positive recurrent.
If $G$ increases,  but    $F(x-1)G(x)\to 0$ as $x\to \infty$, 
 then the CTMC  is  still recurrent. 
If $G$  increases  sufficiently fast so that  $F(x)G(x)\to \infty$ as $x\to \infty$,  then the Markov chain becomes transient and 
can be even  explosive.}
\end{rem}

Before we formulate Theorems~\ref{T2} and~\ref{T3},   we would like to consider 
 an  exponential case,  i.e.    $F(x)=e^{\alpha x}$ and $G(x)=e^{\beta  x}$, where $\alpha, \beta \in \R$.  
Note first that in this case   Theorem~\ref{T1}  yields the following. 
If $\alpha<0$ and $\alpha+\beta<0$,  then  CTMC  $\xi(t)$  is positive recurrent. Also, if 
 either $\alpha<0, \alpha+\beta>0$, or $\alpha>0$, then CTMC $\xi(t)$ is transient.
A direct computation gives  that  the CTMC  is reversible with  the following invariant measure
$e^{\alpha\frac{x(x-1)+y(y-1)}{2}+\beta xy},\,(x,y)\in \Z_{+}^2$, which 
 is summable  if  and only if $\alpha<0, \alpha+\beta<0$.  
Thus,  the sufficient condition of positive recurrence in Theorem~\ref{T1} is also a necessary one in the exponential case.
Note that CTMC $\xi(t)$ in the exponential case   is a particular case of a Markov chain studied  in~\cite{Volkov2015}.
The Markov chain in~\cite{Volkov2015} describes evolution of 
a system of locally interacting  birth-and-death processes labelled by vertices of a finite connected graph.
In terms of~\cite{Volkov2015},   CTMC $\xi(t)$ corresponds to the simplest graph with just two vertices.  
The following proposition is an extract of  results in~\cite{Volkov2015}  complementing   Theorem~\ref{T1} in the exponential case.
\begin{prop}
\label{P0}
\begin{enumerate}
\item[ 1)] If $\alpha<0$ and $\alpha+\beta=0$,  then  CTMC $\xi(t)$ is transient
and does  not  explode.
\item[ 2)]  If either $\alpha>0$, or $\alpha+\beta>0$,  then CTMC  $\xi(t)$ is explosive.
\item[ 3)]  If $\alpha=0$ and $\beta\leq 0$, then both CTMC $\xi(t)$ and DTMC $\zeta(t)$ are null recurrent.
\item[ 4)]  If $\alpha=0$ and $\beta>0$, then DTMC $\zeta(t)$ is transient and CTMC $\xi(t)$  is explosive.
\item[ 5)] Furthermore, 
(i)
 if $\alpha<0$ and $\alpha+\beta\geq 0$,  or, if  $0<\alpha<\beta$, 
then 
$\P(\zeta_1(t)=\zeta_2(t)\ \mbox{infinitely often})=1,$
(ii)  if $\alpha>|\beta|$,  then with probability $1$  eventually a single component of DTMC  $\zeta(t)$ grows while 
 the other component stops changing at all.
\end{enumerate}
\end{prop}

Theorems~\ref{T2} and~\ref{T3} below  are examples of  statements that are similar to Proposition~\ref{P0}.
Namely, these theorems complement Theorem~\ref{T1} 
by providing  more detailed  description of the long term  behaviour  of  the Markov chain under additional 
assumptions  about   functions $F$ and $G$.
Theorem~\ref{T2}  complements Part 1)b) of Theorem~\ref{T1}.
Theorem~\ref{T3} describes a  rather unusual  phenomenon  in a  transient case  specified by  polynomial functions $F$ and $G$. 
\begin{theorem}
\label{T2}
Let functions $F$ and $G$ be positive. Suppose that    function $F$ is non-increasing and $\lim_{x\to \infty}F(x)=0$,
function $G$ is non-decreasing and   $\lim_{x\to \infty}G(x)=\infty$. Suppose also that $\lim_{x\to \infty}F(x)G(x)=\infty$.
Then,
 with probability $1$, 
\begin{itemize}
\item[ 1)]
$\zeta_1(t)=\zeta_2(t)$ for infinitely many $t$;
\item[ 2)] 
if, in addition, 
$\lim_{x\to \infty}\frac{F(x+a)G(x+b)}{F(x)G(x)}=1$  for any  $a, b\in \R$, then given any  $\delta\in (0, 1)$ 
$\zeta(t)\in \{(x,y): \delta x\leq y\leq \delta^{-1}x\}$ 
 for all  but finitely many $t$.
\end{itemize}
\end{theorem}

\begin{theorem}
\label{T3}
Let  $F(x)=(x+1)^{\lambda_1}$ 
 and $G(x)=(x+1)^{-\lambda_2}$, where $0<\lambda_1<\lambda_2$. 
\begin{itemize}
\item[ 1)] 
If $0<\lambda_1\leq 1$, then
 CTMC $\xi(t)$ is transient and non-explosive.
 Further, let 
 $k\in\Z_{+}$ be  such that $\lambda_1+k\lambda_2\leq 1<\lambda_1+(k+1)\lambda_2$. Then,  
 with a  positive probability $\tilde p$ (depending on an initial state),
CTMC  $\xi(t)$ is eventually absorbed by horizontal strip  $\{(x, y)\in \Z_{+}^2: y\leq k\}$
and each of the following  sets $\{t\in \R_{+}: \xi_2(t)= j\}$, $j\leq k$, is unbounded;  
with probability $1-\tilde p$, CTMC $\xi(t)$ is eventually absorbed by vertical strip $\{(x, y): x\leq k\}$ and 
 each of the following  sets $\{t\in \R_{+}: \xi_1(t)= j\}$,  $j\leq k$, is unbounded.

\item[ 2)]
If $\lambda_1>1$, then CTMC $\xi(t)$ is transient and  explodes  with probability $1$. 
Further,  if $\tau_{exp}$ is the time to explosion, then    with probability $1$ 
there exists  a  random integer $m$ and a random time $\tau<\tau_{exp}$  
 such that $\min(\xi_1(t), \xi_2(t))=m$ for all $t\geq \tau$.
In other words, with probability one  there exists  a  random integer $m$
such that the Markov chain explodes by moving eventually along either a horizontal ray $\{(x,y)\in \Z_{+}^2: y= m\}$, or
along a vertical ray $\{(x,y)\in \Z_{+}^2: x= m\}$.
 \end{itemize}
\end{theorem}


\section{Proofs}

\subsection{Proof of Theorem~\ref{T1}}

\paragraph{{\it Proof of Part 1)a) of Theorem~\ref{T1}.}}
There are two cases to consider. If  both  functions $F$ and $G$ are non-increasing  and tend to zero at infinity,
then positive recurrence  of  CTMC $\xi(t)$  is rather obvious and we omit 
 the proof.  
In the second case, where   $\lim_{x\to \infty}F(x)=0$, $\lim_{x\to \infty}G(x)=\infty$  and $\lim_{x\to \infty}F(x-1)G(x)=0$,
we are going to prove positive recurrence of DTMC $\zeta(t)$.
Positive recurrence of the DTMC will yield positive recurrence of CTMC $\xi(t)$ as the transition rates  are uniformly bounded below. 

To prove positive recurrence of the DTMC $\zeta(t)$ we are going to apply 
Theorem 2.2.4 from \cite{FMM} which  is a generalisation 
of the classical  Foster  criterion for  positive recurrence of  irreducible DTMC's  (e.g., Theorem 2.2.3,~\cite{FMM}).
According to this theorem,   DTMC $\zeta(t)$ is positive recurrent,   if   
 there exist  positive   functions $f: \Z_{+}^2\to (0, \infty)$ (the Lyapunov function)  and $\kappa: \Z_{+}^2\to \N=\{1,2,...\}$,  and  $\eps>0$, 
 such that $f(x,y)\to \infty$ as $(x,y)\to \infty$ in any reasonable sense (e.g. $x+y\to \infty$), 
and 
\begin{equation}
\label{Lyap}
m_f(x,y, \kappa(x,y))\leq -\eps\kappa(x,y),
\end{equation}
where $m_f$ is defined by (\ref{mf}),   for all  $(x,y)$
outside a bounded neighbourhood of the origin.
Here we define functions $f$ and $\kappa$  as follows
$$
f(x,y)=\begin{cases}
\alpha x -y, &  0\leq y\leq x,\\
\alpha y-x,& 0\leq x<y,
\end{cases}
$$
where $\alpha>3$,  and 
\begin{equation}
\label{kappa}
\kappa(x,y)=\begin{cases}
1, &  y\neq x,\\
2,& y=x.
\end{cases}
\end{equation}
It is easy to see that $f(x,y)>0$ on $\Z_{+}^2$ and  $f(x,y)\to \infty$ as $x+y\to \infty$.
 Let us verify that inequality  (\ref{Lyap})
is satisfied with these functions.  Without loss of generality,  suppose  that $0\leq y\leq x$.
Notice that, in this case, if $x+y$ is large, then $x$ is also necessarily large (at least $(x+y)/2$).

It is easy to see that if $y<x$,  then inequality   (\ref{Lyap}) becomes
$m_f(x,y, 1)\leq -\eps$, 
or,    equivalently, 
$$(\alpha-\eps)F(x)G(y)-(1+\eps)F(y)G(x)-\alpha+\eps+{\bf 1}_{\{y>0\}}(1+\eps)\leq 0.$$
Monotonicity  of both $F$ and  $G$ imply that  the left side of the preceding display can be bounded by
$$  (\alpha-\eps)F(x)G(x)-(1+\eps)F(0)G(x)-\alpha+1+2\eps,$$
where the first term vanishes and  negative second and third terms dominate for large $x$.

Let us  show that 
\begin{equation}
\label{alpha}
m_f(x,x, 2)\leq -2\eps.
\end{equation}
Starting at $(x,x)$ the Markov chain can reach in two steps the following  states $(x+i, x+j)$, where integers $i$ and $j$ are such $|i|+|j|=2$.
It is easy to see that under  assumptions  of the theorem
 $\lim_{x\to\infty}\gamma(x+a, x+b)=2$,
and $\lim_{x\to\infty}F(x+a)G(x+b)=0$.
This means  that in a finite vicinity of the diagonal located  sufficiently far from the origin  
 the  DTMC  jumps only either down or left with probabilities close to $1/2$, 
and other  jumps can be neglected.
 This yields that 
starting at $(x,x)$, where $x$ is sufficiently large,  
$\zeta(2)$  takes values 
$(x-2, x), (x-1, x-1)$ or $(x, x-2)$ with probabilities converging  
 to  $1/4, 1/2$ and $1/4$ 
respectively, as $x\to \infty$,  and probabilities of other potentially reachable in   two steps states tend to zero in the same limit.
Also, the differences $f(x+i, x+j)-f(x,x)$  are uniformly bounded in $x$. Therefore, 
\begin{align*}
m_f(x, x, 2)&= \frac{f(x-2, x)-2f(x,x)+f(x, x-2)}{4}
+\frac{f(x-1, x-1)-f(x,x)}{2}+C(x)\\
&=\frac{3-\alpha}{2}+C(x),
\end{align*}
where $C(x)\to 0$ as $x\to \infty$, which means that  the left side of~(\ref{alpha}) is less than 
$-2\eps$ for some $\eps>0$ for  all sufficiently large $x$  by the choice of $\alpha$.

\paragraph{{\it Proof of Part 1)b) of Theorem~\ref{T1}.}}
We are going to show transience of DTMC $\zeta(t)$.
Define   $f(x,y)=x+y$  and  $D_a=\{(x,y): x+y\geq a\} \in \Z_{+}^2$, where $a>0$.
Let us  show that  if $a$ is sufficiently large, then there exists $\eps>0$ such that  for all $(x,y)\notin D_a$ 
\begin{equation}
\label{eq1}
m_f(x, y,1)\geq \eps.
\end{equation}
Notice that if $x+y\geq a$ and $0\leq y\leq x$, then necessarily $x\geq a/2$.
It is easy to see that if $0\leq y\leq x$, then   equation~(\ref{eq1}) is equivalent to the following one
$$(F(x)G(y)+F(y)G(x))(1-\eps)-(1+1_{\{y>0\}})(1+\eps)\geq 0,$$
and the left side of the preceding inequality can be bounded below as follows
\begin{align*}
(F(x)G(y)+F(y)G(x))(1-\eps)&-(1+1_{\{y>0\}})(1+\eps) \\
&\geq F(y)G(x)(1-\eps)-2(1+\eps)\\
& \geq F(x)G(x)(1-\eps)-2(1+\eps)\\
&\geq F(a/2)G(a/2)(1-\eps)-2(1+\eps)
\end{align*}
It is easy to see that given $\eps\in (0, 1)$ the right side of the last inequality is positive
for sufficiently large  $a$. 
Thus,  inequality~(\ref{eq1}) holds, which  implies, by Theorem~\ref{227},  that  DTMC $\zeta(t)$ is transient.

\paragraph{{\it Proof of Part 2) of Theorem~\ref{T1}.}}
Recall that  in this part  $F$ is non-decreasing and tends to infinity as $x\to \infty$.
If also   $\lim_{x\to \infty}G(x)=\infty$, then transience  
of the Markov chain is  obvious.
In the rest of the proof we assume that  $G$ is non-increasing and $\lim_{x\to \infty}G(x)=0$.
As in the proof of   Part 1)b), we show  transience of DTMC $\zeta(t)$. 
There are two  cases to consider: $\lim_{x\to \infty}F(x)G(x)=\infty$ and $\lim_{x\to \infty}F(x)G(x)=0$.

Suppose first that $\lim_{x\to \infty}F(x)G(x)=\infty$. 
We are going to show that there exists $\eps>0$ such that  for all $(x,y)\notin D_a$, where $a=a(\eps)$ is sufficiently large, 
inequality~(\ref{eq1}) holds with the same  function $f(x,y)=x+y$ as in the proof of Part 1)b).
Without loss of generality, suppose that  $0\leq y\leq x$, in which case inequality~(\ref{eq1}) is equivalent 
  to the following one
$$(F(x)G(y)+F(y)G(x))(1-\eps)-(1+1_{\{y>0\}})(1+\eps)\geq 0.$$
The left side of the preceding inequality can be bounded below as follows
\begin{align*}
(F(x)G(y)+F(y)G(x))(1-\eps)&-(1+1_{\{y>0\}})(1+\eps)\\
&\geq F(x)G(y)(1-\eps)-2(1+\eps)\\
& \geq F(x)G(x)(1-\eps)-2(1+\eps)\\
&\geq F(a/2)G(a/2)(1-\eps)-2(1+\eps).
\end{align*}
It is easy to see that given $\eps\in (0, 1)$ the right side of the last inequality is positive
for sufficiently large  $a$. Therefore,  by  Theorem \ref{227} DTMC $\zeta(t)$ is transient.

Suppose  now that $ \lim_{x\to \infty}F(x)G(x)=0$. 
 Fix   $\alpha\in (0, 1)$ and  define the following function
\begin{equation}
\label{f2}
f(x, y)=\begin{cases}
\alpha x-y,& 0\leq y< \alpha x,\\
\alpha y-x,& 0\leq x< \alpha y,\\
1, & \mbox{otherwise}.
\end{cases}
\end{equation}
We are going to  show that if  $(x,y)\in A=\{y<\alpha x-C, x\geq a\}\cup \{x< \alpha y - C, y\geq a\}$,
where $C>1$ and $a$ is sufficiently large,  then $m_f(x,y, 1)\geq \eps$ for  $0<\eps<(1-\alpha)/2$.
Due to  symmetry between $x$ and $y$ it suffices to show this bound for $0\leq y<x$,
in which case  inequality  $m_f(x,y,1)\geq \eps$ is equivalent to the following one
$$F(x)G(y)(\alpha-\eps) -F(y)G(x)(1+\eps)+1-\alpha-\eps(1+1_{\{y>0\}})\geq 0.$$ 
If $0< y< x$, then $F(x)G(y)\geq F(x)G(x)$ and $-F(y)G(x)\geq -F(x)G(x)$, therefore 
  the left side of the preceding display   can be bounded below as follows
\begin{align*}
F(x)G(y)(\alpha-\eps)& -F(y)G(x)(1+\eps)+1-\alpha-\eps(1+1_{\{y>0\}})\\
&\geq (\alpha-1-2\eps)F(x)G(x)+1-\alpha-2\eps\\
&\geq (\alpha-1-2\eps)F(a)G(a)+1-\alpha-2\eps,
\end{align*}
and the right side of the last inequality is positive for sufficiently large $a$, as $1-\alpha-2\eps>0$
and  $\lim_{a\to\infty}F(a)G(a)=0$.
Now we apply again  Theorem~\ref{227} with function (\ref{f2}) and set $A$ to finish  the proof.

\subsection{Proof of Theorem~\ref{T2}}

\paragraph{{\it Proof of Part 1) of Theorem~\ref{T2}.}}
Define the following function
$$f(x,y)=
\begin{cases}
x-y,& 0\leq y\leq x,\\
y-x,& y>x.
\end{cases}
$$
If $0\leq y\leq x$, then 
$$
m_f(x,y, 1) =
\frac{F(x)G(y)-F(y)G(x)-1+1_{\{y>0\}}}{\gamma(x,y)}\leq \frac{F(x)G(y)-F(y)G(x)}{\gamma(x,y)}\leq 0,
$$
as $-F(y)\leq -F(x)$ and $G(y)\leq G(x)$. 
 Symmetry between $x$ and $y$ implies that  $m_f(x, y, 1)\leq 0$ holds  in the case $y>x$ as well.
This  yields that 
$\eta(t)=f(\zeta_1(t\wedge \tau), \zeta_2(t\wedge \tau))$, where 
$\tau=\min\{t: \zeta_1(t)=\zeta_2(t)\}$, 
 is a non-negative supermartingale. Therefore, $\eta(t)$   converges almost surely
to a finite limit as $t\to \infty$. This necessarily implies that  $\tau=\min\{t: \zeta_1(t)=\zeta_2(t)\}$ is  almost surely finite 
as $|\eta(t+1)-\eta(t)|=1$ for $t<\tau$, and, hence, with probability $1$
 DTMC $\zeta(t)$  hits the diagonal $y=x$ infinitely 
many times.

\paragraph{{\it Proof of Part 2) of Theorem~\ref{T2}.}}
Given $\delta\in (0, 1)$ 
define  $K_{\delta}=\{(x,y): \delta x\leq y\leq \delta^{-1}x\}$
and   $\sigma=\inf\{t: \zeta(t)\notin K_{\delta}\}$.

\begin{prop}
\label{P11}
There exists $\eps>0$ such that 
$\inf_{(x,y)\in K_{\delta}}\P(\sigma=\infty|\zeta(0)=(x,y))>\eps$.
\end{prop}
{\it Proof of Proposition~\ref{P11}.}
Given $\delta>0$  define     the following functions
$$f(x,y)=
\begin{cases}
y-\delta x,& y\leq x,\\
x-\delta y,& x< y,
\end{cases}
$$
and 
$$\kappa(x,y)=
\begin{cases}
1,& x\neq y,\\
n,& x=y,
\end{cases}
$$
where $n=n(\delta)$ is sufficiently large and to be chosen later.
 We are going to show that 
\begin{equation}
\label{Lyap1}
m_f(x,y,\kappa(x,y))\geq \eps',
\end{equation}
for some $\eps'>0$.
Indeed, if $0<y<x$, then  inequality~(\ref{Lyap1}) becomes  $m_f(x,y,1)\geq \eps'$, which is equivalent to 
$$(1-\eps')F(y)G(x)-(\delta+\eps')F(x)G(y)-1+\delta -2\eps'\geq 0.$$
It is easy to see that the left side of the preceding display can be bounded below as follows
\begin{align*}
(1-\eps')F(y)G(x)&-(\delta+\eps')F(x)G(y)-1+\delta -2\eps'\\
&\geq 
(1-\delta-2\eps')F(x)G(x)-1-2\eps'+\delta\geq 0.
\end{align*}
Due to symmetry between $x$ and $y$  inequality~(\ref{Lyap1})  holds for   $0<x<y$ as well.

If $y=x$  then we  are going to show that, given $0<\delta<1$, there exists 
$n=n(\delta)$ such that 
$m_f(x, x, n)\geq \eps'$,
for some $\eps'>0$. Indeed, assumption  
 $\lim_{x\to\infty}\frac{F(x+a)G(x+b)}{F(x)G(x)}=1$
implies  that 
given integers $n, i$ and $j$ such that $|i|+|j|\leq n$  the  DTMC  jumps  from 
$(x+i, x+j)$  up and  right  with probabilities that tend to $1/2$ as $x\to \infty$.
In turn, this yields that 
starting at $(x,x)$, where $x$ is sufficiently large,  
$\zeta(n)$  takes values 
$(x+k, x+n-k), k=0,\ldots, n$  with probabilities that tend to the  binomial probabilities ${n \choose k}2^{-n}$, $k=0, \ldots, n$
 as $x\to \infty$,  and probabilities of other states reachable in  $n$ steps  tend to zero in the same limit.
Therefore, 
$$m_f(x, x, n)=\E\left(f(x+Y, x+n-Y)\right)-f(x,x)+C(x),$$
where $Y$ is a Binomial random variable with parameters $n$ and $p=1/2$, and $C(x)\to 0$ as $x\to \infty$.
Notice also, that $f(x+a, x+b)=f(x+b, x+a)$ for any $a, b\in \Z$. 
Without loss of generality, assume   that $n=2m+1$. A direct computation (we skip some details) 
gives  that 
\begin{align*}
\E\left(f(x+Y, x+n-Y)\right)-f(x,x)&=\frac{1}{2^{n-1}}\sum\limits_{k=0}^m {n \choose k}(k-\delta(n-k))\\
&=
\frac{1+\delta}{2^{n-1}}\sum\limits_{k=0}^m{n \choose k}k-\frac{\delta n}{2^{n-1}}\sum\limits_{k=0}^m{n \choose k}\\
&=\frac{1+\delta}{2^{n-1}}\left(n2^{n-1}-\frac{n}{2}{2m\choose m}\right)-\delta n\\
&\approx \frac{1-\delta}{2}\left(n-C\frac{1+\delta}{1-\delta}\sqrt{n}\right)>\eps',
\end{align*}
for some $\eps'>0$, if $n$ is large enough.
Given $(x_0,y_0)$ define the following  sequence of random times 
 $n_0=0$ and $n_{t}=n_{t-1}+\kappa(\zeta(n_{t-1}))$, $t\geq 1$, and the following random process $S(t)=f(\zeta(n_t))$, $t\geq 0$.
By construction,  $S(t)\geq 0$ if and only if $\zeta(n_t)\in K_{\delta}$.
Define also $\tau_0=\inf(t: S(t)<0)$. 
It is easy to see that event $\{\tau_0=\infty\}$ implies event $\{\sigma=\infty\}$.
Inequality (\ref{Lyap1}) yields  that  
$\E(S(t+1)-S(t)|S(t))\geq \eps'$ and, therefore,  by Theorem~\ref{219}, we obtain that there exists $\eps>0$ such that 
$\P(\tau_0=\infty|\zeta(0)\in K_{\delta})>\eps$. Consequently, 
$\P(\sigma=\infty|\zeta(0)\in K_{\delta})>\eps$.
Proposition~\ref{P11} is proved. 

Part 1) of the theorem implies that  with probability $1$ DTMC $\zeta(t)$
returns to  set $K_{\delta}$. 
 Define 
$A_m=\{\mbox{the DTMC leaves  set}\, K_{\delta}\, \mbox{at least}\, m \, 
\mbox{times}\}$.
By Proposition~\ref{P11}, we have that  $\P(A_m|A_{m-1})\leq 1-\eps$, where $\eps\in (0, 1)$. Consequently, this 
yields  that  $\P(A_{m})=\P(A_{m}|A_{m-1})\cdots \P(A_1)\leq (1-\eps)^{m}$, so that 
 with probability $1$  DTMC $\xi(t)$ leaves set $K_{\delta}$ finitely many times. The proof of Part 2) of the theorem is finished.

\subsection{Proof of Theorem \ref{T3}}

First we note  that if  $0<\lambda_1<\lambda_2$ then $F(x)=(x+1)^{\lambda_1}\to \infty$, $G(x)=(x+1)^{-\lambda_2}\to 0$
and $F(x)G(x)\to 0$ as $x\to \infty$. 
Therefore transience of the CTMC $\xi(t)$ in both parts of the theorem 
is implied by Theorem~\ref{T1}.

\subsubsection{Proof of Part 1 of Theorem~\ref{T3}}
\label{part1}
The proof is divided on  steps given by 
Propositions~\ref{P1}, \ref{P2} and~\ref{P3}, Corollary~\ref{C1}, and Lemmas~\ref{L1} and~\ref{L2}.
The lemmas  form  the cornerstone of the proof and based on the so called Lyapunov functions approach (e.g., see~\cite{FMM}) widely used 
for study the long term  behaviour  of Markov processes. 

We start with showing non-explosiveness of the CTMC.
\begin{prop}
\label{P1}
Let $F(x)=(x+1)^{\lambda_1}$ and $G(x)=(x+1)^{-\lambda_2}$, where $0<\lambda_1\leq 1$ and $\lambda_2>0$.
Then CTMC $\xi(t)$ is non-explosive with probability $1$.
\end{prop}
{\it Proof of Proposition~\ref{P1}.}
 Let  $\gamma(x,y)$  be a total intensity of jumps of the CTMC at state $(x,y)$.
It is easy to see that
\begin{align*}
\gamma(x,y)&=(x+1)^{\lambda_1}(y+1)^{-\lambda_2}
+(y+1)^{\lambda_1}(x+1)^{-\lambda_2}+
1_{\{x>0\}}+1_{\{y>0\}}\\
&\leq (x+1)^{\lambda_1}
+(y+1)^{\lambda_1}+2\\
&\leq 2(\max(x,y)+1)^{\lambda_1}+2,
\end{align*}
and, hence,
$\gamma^{-1}(x,y)\geq   \left[2(\max(x,y)+1)^{\lambda_1}+2)\right]^{-1}$.
 Let  $(x_n, y_n)$,  $n\in \Z_{+}$, be a  trajectory   of the Markov chain,   such that 
$\lim_{n\to \infty}\max(x_n, y_n)=\infty$,  and consider any of its  subsequences    
$(x_{n_k}, y_{n_k})$, $k\in \Z_{+}$, such that  $\max(x_{n_k}, y_{n_k})=k$.
It is easy to see that 
$$\sum\limits_{n=1}^{\infty}\frac{1}{\gamma(x_n,y_n)}\geq \sum\limits_{k=1}^{\infty}
\frac{1}{\gamma(x_{n_k}, y_{n_k})}\geq  \sum\limits_{k=1}^{\infty}\frac{1}{2((k+1)^{\lambda_1}+1)}=\infty.$$
Thus 
$\sum_{n=1}^{\infty}\gamma^{-1}(x_n,y_n)=\infty$,
 and, hence,  by 
 the well-known criterion of non-explosiveness,    the Markov chain is not  explosive. 
Proposition~\ref{P1} is proved.

\begin{prop}
\label{P2}
Let $F(x)=(x+1)^{\lambda_1}$ and $G(x)=(x+1)^{-\lambda_2}$, where $0<\lambda_1\leq 1$ and $\lambda_1<\lambda_2$.
Let  $\tau_0=\inf\{t: \min(\xi_1(t), \xi_2(t))=0\}$.
Then there exists $\eps>0$ such that 
for any initial state $(x, y)$ 
$$\E(\tau_0|\xi(0)=(x,y))\leq \min(x, y)/\eps.$$
\end{prop}
{\it Proof of Proposition~\ref{P2}.}
Note first that by Proposition~\ref{P1} CTMC $\xi(t)$ is non-explosive.
Denote  $\eta(t)=\min(\xi_1(t), \xi_2(t))$ and define   $Y_t=\eta(t\wedge \tau_0)$.
If  $(\xi_1(t), \xi_2(t))=(x, y)$, where  $0\leq y\leq x$, then $\eta(t)=\xi_2(t)=y$ and  
\begin{align*}
\E(Y(t+dt)-Y(t)|\xi(t)=(x, y))&=
((x+1)^{-\lambda_2}(y+1)^{\lambda_1}-1)dt+\bar{o}(dt)\\
&\leq ((x+1)^{\lambda_1-\lambda_2}-1)dt+\bar{o}(dt)\leq -\eps dt,
\end{align*}
on $\{t<\tau_0\}$, for some $\eps>0$,  and where $\bar{o}(dt)/dt\to 0$ as $dt\to 0$. By the symmetry between $x$ and $y$  we get that 
$$
\E(Y(t+dt)-Y(t)|\xi(t)=(x,y))\leq \left(\left(\max(x,y)+1\right)^{\lambda_1-\lambda_2}-1\right)dt +\bar{o}(dt)\leq -\eps dt,
$$
for  all $(x,y)\in \Z_{+}^2$,  on $\{t<\tau_0\}$.
Proposition~\ref{P2}  is now implied by Theorem~\ref{211} in Appendix.

\smallskip 
  
Proposition~\ref{P1} and Proposition~\ref{P2}  yield the following corollary.
\begin{corollary}
\label{C1}
Under assumptions of Proposition~\ref{P2} set  $\{t\in \R_{+}:\min(\xi_1(t), \xi_2(t))=0\}$ is unbounded with probability $1$.
\end{corollary}
The next lemma states that with a positive probability the Markov chain stays forever in a strip along  one of the coordinate axis.
\begin{lemma}
\label{L1}
Let  $F(x)=(x+1)^{\lambda_1}$ and $G(x)=(x+1)^{-\lambda_2}$, where $0<\lambda_1\leq 1$ and $\lambda_2>0$.
Let  $k\in \Z_{+}$ be such that  $\lambda_1+(k+1)\lambda_2>1$.
Given $N\in \Z_{+}$ define $D_{1, k, N}=\{x\geq N, y\leq k\}$
and $\tau_{1, k, N}=\inf\{t: \xi(t)\notin  D_{1, k, N}\}$. 
Similar, define $D_{2, k, N}=\{x\leq k, y\geq N\}$
and $\tau_{2, k, N}=\inf\{t: \xi(t)\notin  D_{2, k, N}\}$.
If $N$ is sufficiently large then there exists $\delta>0$ such that 
\begin{equation}
\label{dx}
\inf_{(x, y)\in D_{1, k, N}}\P(\tau_{1, k, N}=\infty| \xi(0)=(x,y))>\delta
\end{equation}
and 
\begin{equation}
\label{dy}
\inf_{(x, y)\in D_{2, k, N}}\P(\tau_{2, k, N}=\infty| \xi(0)=(x, y))>\delta.
\end{equation}
\end{lemma}
Lemma~\ref{L1} is proved in Section~\ref{ProofL1}.

We are interested in the minimal $k$ satisfying the requirement of Lemma~\ref{L1}.
Namely, let  $k_{min}$ be such that 
$\lambda_1+\lambda_2k_{min}\leq 1< \lambda_1+\lambda_2(k_{min}+1)$.
As the Markov chain is transient, we can assume for the rest of the proof that 
$N$ is so large that  i) sets $D_{1, k_{min}, N}$ and $D_{2, k_{min}, N}$
are disjoint; ii)  bounds (\ref{dx}) and (\ref{dy}) hold. 

\begin{prop}
\label{P3}
With a positive probability $\tilde p$, depending on $\xi(0)$,
 CTMC $\xi(t)$ is  eventually  absorbed
 by horizontal strip $D_{1, k_{min}, N}$, and with probability $1-\tilde p$ CTMC $\xi(t)$
is eventually absorbed by  vertical strip $D_{2, k_{min}, N}$.
\end{prop}
{\it Proof of Proposition~\ref{P3}.}
Note first that by Corollary~\ref{C1}      CTMC $\xi(t)$
returns to  set  $\{x\geq N, y\leq k_{min}\}\cup \{x\leq k_{min}, y\geq N\}$ with probability $1$.
Further, by Lemma~\ref{L1}, if the Markov chain is in either of these strips, then it remains there  with a probability bounded away from zero.
Consequently, with probability $1$  CTMC $\xi(t)$  is eventually absorbed by the union of these  strips.
This can be shown in the same way as the similar fact  in   the proof of  Part 2) of Theorem~\ref{T2} (i.e. absorption by cone $K_{\delta}$).
Finally, it is obvious that absorption by strip   $\{x\geq N, y\leq k_{min}\}$ and absorption by strip  $\{x\leq k_{min}, y\geq N\}$ 
are mutually exclusive events, as the strips are disjoint by assumption.
Proposition~\ref{P3}  is proved.

\begin{lemma}
\label{L2}
Define  $\tau_{k, 1}=\inf(t\geq 0: \xi_1(t)=k)$ and $\tau_{k, 2}=\inf(t\geq 0: \xi_2(t)=k)$.
If  $0<\lambda_1< 1, \lambda_2>0$ and integer $k\geq 1$ are such that  $\lambda_1+k\lambda_2\leq 1$,  then
$$\P(\tau_{k, 1}<\infty|\xi_1(0)=0)=\P(\tau_{k, 2}<\infty|\xi_2(0)=0)=1.$$
\end{lemma}
 Lemma~\ref{L2} is proved  in Section~\ref{ProofL2}. Now we  use this lemma  to finish the proof.
 Lemma~\ref{L2} and  Corollary~\ref{C1} 
yield that  if  CTMC $\xi(t)$  is absorbed by horizontal strip  $\{x\geq N, y\leq k_{min}\}$, then 
 it visits each of the following sets $y\equiv i$,  $i=0,\ldots, k_{min}$, infinitely many times.
Similar, if CTMC $\xi(t)$  is absorbed by vertical strip
 $\{x\leq k_{min}, y\geq N\}$,  it visits each of the following sets $x\equiv i$,  $i=0,\ldots, k_{min}$, infinitely many times.

Part 1) of Theorem~\ref{T3} is now proved.

\subsubsection{Proof of Part 2) of Theorem~\ref{T3}}

Given  $m\in \Z_{+}$  and  $0<\nu<\lambda_1-1$,  define the following function
\begin{equation}
\label{f}
f(x, y)=\begin{cases}
x^{-\nu},& y=m,\, x>0\\
1,& y\neq m \,\, \text{or} \,\, x=0.
\end{cases}
\end{equation}
It is easy to see that 
\begin{align}
\nonumber
\G f(x, m)&=\left(\frac{1}{(x+1)^{\nu}}-\frac{1}{x^{\nu}}\right)\frac{(x+1)^{\lambda_1}}{(m+1)^{\lambda_2}}+
\left(\frac{1}{(x-1)^{\nu}}-\frac{1}{x^{\nu}}\right)\\
&+
\left(1-\frac{1}{x^{\nu}}\right)\left(\frac{(m+1)^{\lambda_1}}{(x+1)^{\lambda_2}}+1\right), \label{eps}\\
&\leq -C_1x^{-\nu-1+\lambda_1}+\nu x^{-\nu-1}+\frac{(m+1)^{\lambda_1}}{(x+1)^{\lambda_2}}+1\leq -\eps,\nonumber
\end{align}
for some $\eps>0$ and for all  $x\geq N_m$, where $N_m$ is sufficiently large.
Bound~(\ref{eps})  implies that 
conditioned to stay in set  $K_{m, N}$   CTMC $\xi(t)$ explodes, with  a positive  probability  depending on $m$,
 by Theorem  1.12,~\cite{MenPetr}. 
By symmetry between $x$ and $y$ we immediately obtain the same for any vertical ray $\{y\geq N_m, x=m\}$.
Let  $\tau_{exp}$ be  the time to explosion, 
$\tau_0=\inf\{t: \min(\xi_1(t), \xi_2(t))=0\}$ (as in Proposition~\ref{P2})
and $\tau=\min(\tau_{exp}, \tau_0)$.
One can show, by repeating verbatim the proof of Proposition~\ref{P2}, that there exists $\eps>0$ such that 
$\E(\tau|\xi(0)=(x,y))\leq \min(x,y)/\eps$. 
This bound and 
conditional explosion along a horizontal and a vertical ray yield  that  $\P(\tau_{exp}<\infty)=1$.
Next, it is easy to see  that $\min(\xi_1(t), \xi_2(t))$ jumps with uniformly bounded rates, therefore 
it changes finitely many times before explosion.
This yields that the Markov chain eventually explodes  being absorbed by either a  horizontal ray $\{y=const\}$ or 
a vertical ray $\{x=const\}$.



\subsection{Proof of Lemma \ref{L1}}
\label{ProofL1}

Due to symmetry between $x$ and $y$ it suffices to prove bound~(\ref{dx}) only.
It should be noted that the proof is reminiscent of the proof of the well  known criteria for transience  of a countable Markov 
chain (e.g., Theorem 2.2.2, \cite{FMM}).  In particular, it consists in constructing a {\it bounded} positive function $f$ 
such that random process $f(\xi(t))$ is supermartingale.

Fix  an integer $k\geq 1$ such that 
$0<\lambda_1\leq 1<\lambda_1+(k+1)\lambda_2$.
Suppose there exists  a positive function $f_k$ on  $\Z_{+}^2$ such that  
\begin{enumerate}
\item  $\max_{y\leq k}f_k(x, y)\to 0$ as $x\to \infty$,
\item $\G f_k(x, y)\leq 0$ for all $(x, y)\in \{x\geq N, y\leq k\}$, 
\item $\sup_{x\geq N}\max_{i=0,..., k}f_k(x, i)\leq N^{-\beta}$,  where 
$N, \beta>0$, and 
\item $f_k(x,y)=1$ for all $(x,y)\notin \{x\geq N, y\leq k\}$.
\end{enumerate}
Define $\tau=\inf(t: \xi(t)\notin \{x\geq N, y\leq k\})$.
The properties of $f_k$ imply   that random process 
$\eta(t)=f_k(\xi_1(t\wedge \tau), \xi_2(t\wedge \tau))$
is a positive supermartingale and, hence, it almost surely 
converges  to a finite limit $\eta_{\infty}$ that can  take only   
 values $1$ and $0$. 
By Fatou's Lemma 
\begin{align*}
\E(\eta_{\infty}|\xi(0)=(x, y))=\P(\tau<\infty| \xi(0)=(x, y))&
\leq \E(\eta(0)| \xi(0)=(x, y))\\
&=f_k(x,y)\leq N^{-\beta},
\end{align*}
for all   $(x, y)\in \{x\geq N, y\leq k\}$, and, hence, 
$\P(\tau=\infty|\xi(0)=(x, y))\geq 1- N^{-\beta}$, 
 for all $(x,y)\in \{x\geq N, y\leq k\}$.
In the rest of the proof we provide   functions $f_k$.

\paragraph{{\it Function   $f_0$.}}
Fix $0<\nu<\lambda_1+\lambda_2-1$ and define the following function  
\begin{equation}
\label{f0}
f_0(x, y)=\begin{cases}
1,& y>0 \,\, \text{or} \,\, x=0,\\
x^{-\nu},& y=0,\, x>0.
\end{cases}
\end{equation}
The following bound is obvious 
\begin{equation}
\label{b0}
\sup\limits_{x\geq N}f_0(x, 0)\leq N^{-\nu}.
\end{equation}
Let us show that, if  $x\geq N$, where $N$ is sufficiently large, then
$\G f_0(x,0)\leq 0$.
Indeed, a direct computation gives that 
\begin{align*}
\G f_0(x, 0)&=\left(\frac{1}{(x+1)^{\nu}}-\frac{1}{x^{\nu}}\right)(x+1)^{\lambda_1}
+\left(\frac{1}{(x-1)^{\nu}}-\frac{1}{x^{\nu}}\right)+
\left(1-\frac{1}{x^{\nu}}\right)(x+1)^{-\lambda_2}\\
&\leq -C_1x^{-\nu-1+\lambda_1}+\nu x^{-\nu-1}+C_2x^{-\lambda_2}\leq 0,
\end{align*}
for all sufficiently large $x$, as $\lambda_1+\lambda_2-1>\nu$.
\paragraph{{\it Functions $f_k$, $k\geq 1$.}}
If $k=1$, then we  define 
\begin{equation}
\label{f1}
f_1(x,y)=\begin{cases}
1,&y\geq 2\,\,\text{or}\,\,  x=0,\\
x^{-\nu_1}, & y=1, x>0,\\
x^{-\nu_1}-x^{-\nu_1-\nu_2}, & y= 0,\,  x>0,
\end{cases}
\end{equation}
where $\nu_1>0$ and $\nu_2>0$ are such that $\nu_1+\nu_2<\lambda_2$ and $\lambda_1+\lambda_2+\nu_2>1$
(it is easy to see that such numbers $\nu_1$ and $\nu_2$ exist).
If  $k\geq 2$, then we define 
\begin{equation}
\label{fk}
f_k(x, y)=\begin{cases}
1,&y\geq k+1\,\,\text{or}\,\,  x=0,\\
x^{-\nu_1}, & y=k, x>0,\\
x^{-\nu_1}-x^{-\nu_1-\nu_2}, & y= k-1,\,  x>0,\\
x^{-\nu_1}-x^{-\nu_1-\nu_2}-\ldots-x^{-\nu_1-\nu_2-\ldots-\nu_{k+1-y}},& y=0,\ldots, k-2, \,  x>0,
\end{cases}
\end{equation}
where positive real numbers $\nu_1, \ldots, \nu_{k+1}$ satisfy the following system of inequalities
\begin{equation}
\label{Cond20}
\begin{cases}
\lambda_2>\nu_1+\nu_2,&\\
0<\nu_{i}<\lambda_2, \,  i=3,\ldots, k+1,&\\
1<\lambda_1+\lambda_2+\nu_2+\ldots+\nu_{k+1}.
\end{cases}
\end{equation}
It is easy to see that  system of  inequalities~(\ref{Cond20}) has  many solutions 
and  for all $k\geq 1$ the following bound holds
\begin{equation}
\label{bk}
\sup\limits_{x\geq N}\max\limits_{0\leq i \leq k}f_k(x, i)\leq \sup\limits_{x\geq N}f_k(x, k)\leq N^{-\nu_1}.
\end{equation}
A direct computation gives that  
\begin{align*}
\G f_1(x, 1)&=\left(\frac{1}{(x+1)^{\nu_1}}-\frac{1}{x^{\nu_1}}\right)\frac{(x+1)^{\lambda_1}}{2^{\lambda_2}}
+\left(\frac{1}{(x-1)^{\nu_1}}-\frac{1}{x^{\nu_1}}\right) \\
&+\left(1-\frac{1}{x^{\nu_1}}\right)\frac{2^{\lambda_1}}{(x+1)^{\lambda_2}}-\frac{1}{x^{\nu_1+\nu_2}},
\end{align*}
and, hence,  
$$
\G f_1(x, 1)\leq  -C_1x^{-1-\nu_1+\lambda_1}
+C_2x^{-\lambda_2}-x^{-\nu_1-\nu_2}\leq 0, 
$$
$$\G f_1(x, 0)\leq  x^{\lambda_1}(-C_1x^{-1-\nu_1}+C_2x^{-\nu_1-\nu_{2}-\lambda_1-\lambda_2})\leq 0,$$ 
 for all sufficiently large 
$x$,  as $\lambda_2>\nu_1+\nu_2$ and $\lambda_1+\lambda_2+\nu_2>1$. 

If $k\geq 2$, then 
a direct computation gives that  
\begin{align*}
\G f_k(x, k)&\leq \left(\frac{1}{(x+1)^{\nu_1}}-\frac{1}{x^{\nu_1}}\right)\frac{(x+1)^{\lambda_1}}{(k+1)^{\lambda_2}}
+\left(\frac{1}{(x-1)^{\nu_1}}-\frac{1}{x^{\nu_1}}\right)\\
&+\left(1-\frac{1}{x^{\nu_1}}\right)\frac{(k+1)^{\lambda_1}}{(x+1)^{\lambda_2}}-\frac{1}{x^{\nu_1+\nu_2}}\\
&\leq  -C_1x^{-1-\nu_1+\lambda_1}
+C_2x^{-\lambda_2}-x^{-\nu_1-\nu_2}\leq 0,
\end{align*}
for sufficiently large $x$, 
as $\lambda_2>\nu_1+\nu_2$. Further, 
given $i=2,\ldots, k$, we get in a similar way that 
$$
\G f_k(x, k+1-i)\leq  x^{\lambda_1}(-C_1x^{-1-\nu_1}
+C_2x^{-\nu_1-\ldots-\nu_i-\lambda_1-\lambda_2}-x^{-\nu_1-\nu_2-\ldots-\nu_{i+1}-\lambda_1}).
$$
Notice that the second inequality  of~(\ref{Cond20}) implies that 
$$-\nu_1-\ldots-\nu_i-\lambda_1-\lambda_2<-\nu_1-\nu_2-\ldots-\nu_{i+1}-\lambda_1,$$
and, hence,  
$\G f_k(x, k-i+1)\leq 0$, provided that $x$ is sufficiently large.

Finally, the bottom inequality in~(\ref{Cond20}) implies that
$$
\G f_k(x, 0)\leq  x^{\lambda_1}(-C_1x^{-1-\nu_1}
+C_2x^{-\nu_1-\ldots-\nu_{k+1}-\lambda_1-\lambda_2})\leq 0,$$
  for sufficiently large $x$.

The lemma is proved.

\subsection{Proof of Lemma~\ref{L2}}
\label{ProofL2}
Due to  symmetry between $x$ and $y$ it  suffices to prove only that  $\P(\tau_{k, 2}<\infty|\xi_2(0)=0)=1$.
It should be noted that the proof  is reminiscent of the proof of the well-known criteria for  recurrence of a countable Markov chain
(e.g., Theorem 2.2.1, \cite{FMM}). In particular, it 
consists in constructing an {\it unbounded} positive function $g$ such that random process $g(\xi(t))$  is a supermartingale.

Given an integer $k\geq 1$,
we are going to construct function  $g_k$ satisfying the following conditions
\begin{enumerate}
\item   $g_k(x,i)\to \infty$ as $x\to \infty$  for all  $0\leq i< k$,
\item $g_k(x,y)=1$ on $\{x=0\}\cup \{y\geq k\}$, 
\item $\G g_k(x, y)\leq 0$ for all $(x,y)\in \{x\geq N, y\leq  k-1\}$, where   $N>0$.
\end{enumerate}
Properties of such function $g_k$  imply that the random process $\eta_k(t)=g_k(\xi(t\wedge \tau_{k, 2}))$ 
 is a positive supermartingale and, hence, converges  almost surely. 
If $(x,y)\in \{x\geq N, y\leq  k-1\}$, then the Markov chain jumps to the right with a rate that is  approximately equal to $x^{\lambda_1}$
for sufficiently large  $x\to \infty$, while
rates of  jumps down, up or left are uniformly bounded over states $(x,y)\in \{x\geq N, y\leq  k-1\}$.
It  means that conditioned to stay in strip  $\{x\geq N, y\leq  k-1\}$  component $\xi_1(t)$  tends to infinity as $t\to \infty$
and, by construction, so does  $\eta_k(t)$, which 
 contradicts  its  convergence, unless  $\P(\tau_{k, 2}<\infty)=1$.

In the rest of the proof we construct the functions $g_k$, $k\geq 1$. Note that  in what follows we write
$\psi(x)\approx \phi(x)$ for all  sufficiently large $x$, if $\lim_{x\to \infty}\psi(x)/\phi(x)=1$.

\paragraph{{\it Function $g_1$.}} Suppose that   $\lambda_1+\lambda_2\leq 1$ and define 
$$g_1(x,y)=\begin{cases}
1,& y\geq 1\,\,\text{or}\,\,  x=0,\\
x^{\nu_1}, & x>0, y=0,
\end{cases}
$$
where  $0<\nu_1<1$.
It is easy to see that 
$$\G g_1(x,0)\approx \nu_1x^{\nu_1-1+\lambda_1}-x^{\nu_1-\lambda_2},$$
for all sufficiently large $x$.
If $\lambda_1+\lambda_2<1$, then   $\nu_1-1+\lambda_1<\nu_1-\lambda_2$, hence, 
 $\G g_1(x, 0)\leq 0$.  If $\lambda_1+\lambda_2=1$, then 
$\G g_1(x, 0)\approx  (\nu_1-1)x^{\nu_1-\lambda_2}<0$ for all  sufficiently large $x$, as $\nu_1<1$.

\paragraph{{\it Function  $g_2$.}} 
If $\lambda_1+2\lambda_2\leq 1$,  then we define 
$$
g_2(x,y)=\begin{cases}
1,&y\geq 2\,\,\text{or}\,\,  x=0,\\
x^{\nu_1}, & y=1, x>0,\\
x^{\nu_1}+B_1x^{\nu_1-\nu_2},& y=0, x>0,
\end{cases}
$$
where 
\begin{equation}
\label{k=2}
\begin{cases}
\lambda_2<\nu_2<\nu_1,\,  \nu_2<1-\lambda_1-\lambda_2,\,   B_1=1,& \text{if}\,  \lambda_1+2\lambda_2<1, \\
\nu_2=\lambda_2<\nu_1<B_1<2^{\lambda_1},& \text{if}\,  \lambda_1+2\lambda_2=1. 
\end{cases}
\end{equation}
It is easy to see that 
$\G g_2(x, 0)\approx \nu_1x^{\nu_1-1+\lambda_1}-B_1x^{\nu_1-\nu_2-\lambda_2}$ for all  sufficiently large $x$.
If $\lambda_1+2\lambda_2=1$, then $\G g_2(x, 0)\leq 0$, 
because of the bottom line in condition~(\ref{k=2}).
If $\lambda_1+2\lambda_2<1$, then the upper line in condition~(\ref{k=2}) yields that 
 $\nu_1-1+\lambda_1<\nu_1-\nu_2-\lambda_2$, and, hence,  $\G g_2(x,0)\leq 0$ for all sufficiently large $x$. 

Further, it is easy to see that 
$\G g_2(x, 1)\approx 2^{-\lambda_2}
\nu_1x^{\nu_1-1+\lambda_1}+B_1x^{\nu_1-\nu_2} -2^{\lambda_1}x^{\nu_1-\lambda_2}$ for all  sufficiently large $x$.
If $\lambda_1+2\lambda_2=1$, then both positive terms  are  smaller than $2^{\lambda_1}x^{\nu_1-\lambda_2}$, as 
$\nu_1-1+\lambda_1<\nu_1-\lambda_2$ and $B_1<2^{\lambda_1}$ respectively.
If $\lambda_1+2\lambda_2<1$, then the negative term dominates both positive terms because 
$\nu_1-1+\lambda_1<\nu_1-\lambda_2$ (as $\lambda_1+\lambda_2<1$), 
and $\nu_1-\nu_2<\nu_1-\lambda_2$ (as $\lambda_2<\nu_2$). Hence, we have again that $\G g_2(x, 1)\leq 0$ for all sufficiently large $x$.

\paragraph{{\it Functions  $g_k, k\geq 3$.}}
If $\lambda_1+k\lambda_2\leq 1$, where  $k\geq 3$, then we  define function $g_k$ as follows 
$$
g_k(x,y)=\begin{cases}
1,&y\geq k\,\,\text{or}\,\,  x=0,\\
x^{\nu_1}, & y=k-1, x>0,\\
x^{\nu_1}+B_1x^{\nu_1-\nu_2},& y=k-2, x>0,\\
x^{\nu_1}+B_1x^{\nu_1-\nu_2}+\ldots+B_{i-1}x^{\nu_1-\nu_2-\ldots-\nu_{i}},& y=k-i,\, i=3,\ldots,k, \,  x>0,
\end{cases}
$$
where 
\begin{equation}
\label{Cond1}
\begin{cases}
B_i=1,\, i=1,\ldots, k-1,\, \lambda_2<\nu_{i},  i=2, \ldots, k, &\\
 \nu_{2}+\ldots+\nu_k<\min(1-\lambda_1-\lambda_2, \nu_1),
& {\rm if}\,\,   \lambda_1+\lambda_2k<1, 
\end{cases}
\end{equation}
and
\begin{equation}
\label{Cond2}
\begin{cases}
B_1<k^{\lambda_1}, \lambda_2(k-1)<\nu_1<B_{k-1},&\\
  B_i<(k-i+1)^{\lambda_1}B_{i-1}, \, i=2,\ldots, k-2, 
\nu_i=\lambda_2,\, i=2,\ldots, k,
& {\rm if}\,\,  \lambda_1+\lambda_2k=1.
\end{cases}
\end{equation}
A direct computation gives that 
\begin{align*}
\G g_k(x,0)&\approx \nu_1x^{\nu_1-1+\lambda_1}
-B_{k-1}x^{\nu_1-\nu_2-\ldots-\nu_k-\lambda_2}\\
&=x^{\nu_1-1+\lambda_1}(\nu_1-B_{k-1}x^{1-\lambda_1-\lambda_2-\nu_2-...-\nu_k}),
\end{align*}
for all sufficiently large $x$, where 
if $\lambda_1+\lambda_2k=1$, then the right hand side is  $x^{\nu_1-1+\lambda_1}(\nu_1-B_{k-1})<0$
by condition (\ref{Cond2}), and if  $\lambda_1+\lambda_2k<1$, then the right hand side is negative 
 by condition~(\ref{Cond1}).
 
Further, a direct computation gives that 
\begin{align*}
\G g_k(x,k-i)
&\approx \nu_1(k-i+1)^{-\lambda_2}x^{\nu_1-1+\lambda_1}-(k-i+1)^{\lambda_1}B_{i-1}
x^{\nu_1-\nu_2-\ldots-\nu_i-\lambda_2}\\
&+B_ix^{\nu_1-\nu_2-...-\nu_i-\nu_{i+1}},
\end{align*}
for   $i=2,\ldots,k-2$, for all sufficiently large $x$.
As before,  consider two cases. If $\lambda_1+\lambda_2k=1$, then 
$$x^{\nu_1-1+\lambda_1}(\nu_1-(k-i+1)^{\lambda_1}B_{i-1}
x^{1-\lambda_1-\lambda_2-\nu_2-\ldots-\nu_i})<0,$$
as $1-\lambda_1-\lambda_2-\nu_2-\ldots-\nu_i=1-\lambda_1-i\lambda_2>0$, so that the first positive term is asymptotically  dominated by the negative one. 
Also, comparing the negative term with the second positive one we get that 
\begin{align*}
B_ix^{\nu_1-\nu_2-...-\nu_i-\nu_{i+1}}&-(k-i+1)^{\lambda_1}B_{i-1}
x^{\nu_1-\nu_2-\ldots-\nu_i-\lambda_2}\\
&=x^{\nu_1-i\lambda_2}(B_i-(k-i+1)^{\lambda_1}B_{i-1})<0,
\end{align*}
by~(\ref{Cond2}).
If $\lambda_1+\lambda_2k<1$, then 
condition (\ref{Cond1}) implies that  
$\nu_1-1<\nu_1-\nu_2-\ldots-\nu_i-\lambda_1-\lambda_2$ and
$\nu_1-\nu_2-\ldots-\nu_i-\lambda_1-\lambda_2> \nu_1-\nu_2-\ldots-\nu_i-\nu_{i+1}-\lambda_1$,
so that  $\G f_{k}(x, k-i)\leq 0$ for all sufficiently large $x$. 

Finally, we get  that 
$$
\G g_k(x,k-1) \approx \nu_1k^{-\lambda_2}x^{\nu_1-1+\lambda_1}+B_1x^{\nu_1-\nu_2}-k^{\lambda_1}x^{\nu_1-\lambda_2}\leq 0,
$$
for all sufficiently large $x$.
Indeed, if $\lambda_1+\lambda_2k\leq 1$, $k\geq 3$,  then $\nu_1-1+\lambda_1<\nu_1-\lambda_2$,
 so that   term $k^{\lambda_1}x^{\nu_1-\lambda_2}$ is  larger (for sufficiently large $x$) 
than $\nu_1x^{\nu_1-1+\lambda_1}$. 
To deal with another positive term in the preceding display, we  consider two cases. 
If $\lambda_1+\lambda_2k<1$, then $\nu_1-\nu_2<\nu_1-\lambda_2$,
because of condition~(\ref{Cond1}). If $\lambda_1+\lambda_2k=1$,  then $\nu_1-\nu_2=\nu_1-\lambda_2$, but $B_1<k^{\lambda_1}$.
Thus, in both cases $B_1x^{\nu_1-\nu_2}<k^{\lambda_1}x^{\nu_1-\lambda_2}$ for all sufficiently large $x$.

The lemma is proved.

\section*{Appendix}

\begin{rem}
\label{R5}
{\rm 
It should be noted that our results imply transience  of CTMC $\xi(t)$ 
in the case of polynomial functions  $F(x)=(x+1)^{\lambda_1}$ and $G(x)=(x+1)^{-\lambda_2}$ for {\it any}  $\lambda_1, \lambda_2>0$.
Indeed, if $0<\lambda_1<\lambda_2$ then, as it is mentioned at the beginning of the proof of Theorem~\ref{T3},
Theorem~\ref{T1} applies. 
If  $0<\lambda_1\leq 1$ and $\lambda_2>0$,  then  Lemma~\ref{L1} implies transience.
If $\lambda_1>1$, then (whatever $\lambda_2$ is) transience is implied the criteria for transience of a countable Markov chain
(e.g. Theorem 2.2.2, \cite{FMM})  which applies  in this case  with Lyapunov function~(\ref{f}).
Further, condition 
 $0<\lambda_1<\lambda_2$ in Theorem~\ref{T3} is not necessary to show {\it just} transience. We essentially  use 
this  condition in both parts of  Theorem~\ref{T3} to  describe  how {\it exactly}  the Markov chain   escapes  to infinity. 
}
\end{rem}
For the reader's convenience we provide some facts that were 
 used in our paper. Theorem \ref{227} is a version of Theorem 2.2.7, \cite{FMM}, Theorem \ref{219} 
is a version of Theorem 2.1.9, \cite{FMM}, and
Theorem~\ref{211}  is Lemma 7.3.6 in \cite{MPW}.
 
\begin{theorem}
\label{227}
{\rm (Theorem 2.2.7, \cite{FMM})}.
Let $\eta(t)$ be an irreducible aperiodic  discrete time  Markov chain on a countable  space ${\cal A}$.
For $\eta(t)$ to be transient, it suffices that there exist
a positive function $f(\eta)$, $\eta\in {\cal A}$,  a bounded  positive integer valued function $\kappa(\eta)$,
 $\eta\in {\cal A}$, and numbers $\eps, C>0$ such
that, setting  $A_C=\{\eta\in {\cal A}: f(\eta)\geq C\}\neq \emptyset$, the following conditions hold:

1) $\sup_{\eta\in {\cal A}} \kappa(\eta)<\infty$;

2)  $\E(f(\eta(t+\kappa(\eta))|\eta(t)=\eta)-f(\eta)\geq \eps$
for all  $\eta\in A_C$;

3) for some $d>0$, the inequality $|f(\eta')-f(\eta'')|>d$ implies that the transition probability  from $\eta'$ to $\eta''$ is zero.
\end{theorem}

\begin{theorem}
\label{219}
{\rm (Theorem 2.1.9, \cite{FMM}).}
Let $\eta(t),\, t\in \Z_{+},$ be 
$\R_{+}$-valued  process adapted to a filtration $({\cal F}_t, t\in \Z_{+})$.
Define  $\tau_C=\min(t\geq 1: \eta(t)\leq C\}$, where $C>0$.  Suppose that its  jumps $\eta(t+1)-\eta(t),\, t\in \Z_{+}$, are uniformly bounded
and there exists $\eps>0$  such that  
$\E(\eta(t+1)|{\cal F}_{t})\geq \eta(t)+\eps$, on $\{t\leq \tau_C\}$, and $\eta(0)>C$.
Then $\P(\tau_C=\infty)>0$.
\end{theorem}

\begin{theorem}
\label{211}
 {\rm (Lemma  7.3.6, \cite{MPW}).}
 Let $(\eta(t), t\in \R_{+})$ be an $\R_{+}$-valued   process adapted to a filtration $({\cal F}_t, t\in \R_{+})$
and let $\tau=\inf(t:\eta(t)=0)$. Suppose that there exists $\eps>0$  such that
$\E(\eta(t+dt)-\eta(t)| {\cal F}_{t-})\leq -\eps dt$, on $\{t\leq \tau\}$. Then $\E(\tau|{\cal F}_0) \leq \eps^{-1}\eta(0)$.
\end{theorem}


\begin{thebibliography}{1}
\bibitem{Anderson} Anderson W. (1991). Continuous time Markov chains: an application oriented approach. Springer Verlag. 
\bibitem{Klebaner}  Barbour, A.D., Hamza, K., Kaspi, H.,  and Klebaner, F.C. (2015). 
Escape from the boundary in Markov population processes. {\it Advances in Applied Probability}, 
{\bf 47}, 4, pp.~1190--1211.
\bibitem{Becker}
Becker, N. G. (1970). A stochastic model for two interacting populations.
{\it Journal of Applied Probability}, {\bf 7}, pp.~544-564.
\bibitem{Chung}  Chung, K. L. (1967). {\it Markov Chains with Stationary Transition Probabilities}. 
2nd ed., in: Die Grundlehren der Mathematischen Wissenschaften Band 104, Springer-Verlag New York, Inc., New York. 
\bibitem{FMM} Fayolle, G., Malyshev, V., and Menshikov, M. (1995). 
{\it Constructive topics in the theory of countable Markov chains}. Cambridge University Press.
\bibitem{Feller} Feller, W. (1968).  {\it An Introduction to Probability Theory and its Applications}. Volume 1, 3rd Edition. John Wiley$\&$Sons, Inc.
\bibitem{Gauss} Gauss,  G. F., Smagardova, N.P., and Witt, A.A. (1936).  Further studies of interaction between 
predator and prey. {\it Journal of Animal Ecology},  {\bf 5}, pp.~1--18.  
\bibitem{Karlin0} Karlin, S., and Taylor, H.M.  (1975). {\it A First Course in Stochastic Processes}. 2nd Edition, Elsevier. 
\bibitem{Kolm72} Kolmogorov, A. N. (1972). The quantitative measurements of mathematical models 
in the dynamics of populations. {\it Problems in Cybernetics}, {\bf 25}, pp.~100--106. (In Russian).
\bibitem{Liggett} Liggett, T. (2010). {\it Continuous time Markov processes -- an Introduction}.  Graduate Studies in Mathematics, American Mathematical Society.
\bibitem{MPW} Menshikov, M.V., Popov, S.  and Wade, A.R. (2016). Non-homogeneous Random Walks:
Lyapunov Function Methods for Near-Critical Stochastic Systems. Cambridge University press.
\bibitem{MenPetr} Menshikov, M. and Petritis (2014).
 Explosion, implosion, and moments of passage times for continuous-time Markov chains: semimartingale approach.
{\it Stochastic Processes and Their Applications,}  {\bf 124},  pp.~2388--2414.
\bibitem{Reuter}
Reuter, G. E. H. (1961).  Competition processes. In: Neyman J. (Ed.)
 Proceedings of The Fourth Berkeley Symposium on Mathematical Statistics and Probability, v.II: Contributions to Probability Theory. 
University of California Press, Berkeley. 
\bibitem{Ridler} Ridler-Rowe, C.J. (1978). On competition between two species. {\it Journal of Applied Probability},
{\bf 15}, pp.~457--465. 


\bibitem{Sigmund}  Sigmund, L. (2007). 
Kolmogorov and population dynamics. In: E. Charpentier, A. Lesne, and N. Nikolski (Eds.) Kolmogorov’s Heritage in Mathematics. Springer Verlag. 

\bibitem{Volkov2015} Shcherbakov, V. and Volkov, S. (2015).
Long term behaviour  of locally interacting birth-and-death processes.
 {\it  Journal of Statistical Physics,} ~{\bf 158}, N1, pp.~132--157.
\end{thebibliography}
\end{document}